\documentclass{article}
\usepackage{amsmath}
\usepackage{amssymb}
\usepackage{a4}
\usepackage{theorem}
\newcommand{\nx}{\negthinspace}
\newcommand{\bzero}{\mathbf{0}}
\newcommand{\mud}{\mu\left(d\right)}
\newcommand{\CC}{\mathbb{C}}

\newcommand{\FF}{\mathbb{F}_2}
\newcommand{\OO}{\mathcal{O}}

\newcommand{\PP}{\mathbb{P}}
\newcommand{\Pthree}{\mathbb{P}_3}
\newcommand{\tPthree}{\tilde{\mathbb{P}}_3}
\newcommand{\PthreeC}{\Pthree\left(\CC\right)}
\newcommand{\tPP}{{\tilde{\mathbb{P}}}}
\newcommand{\ZZ}{\mathbb{Z}}
\newcommand{\NN}{\mathbb{N}}
\newcommand{\tD}{\tilde{D}}
\newcommand{\oC}{\overline{C}}
\newcommand{\oD}{\overline{D}}
\newcommand{\oQ}{\overline{Q}}

\newcommand{\ose}{\overline{e}}
\newcommand{\og}{\overline{g}}
\newcommand{\tH}{\tilde{H}}
\newcommand{\tS}{\tilde{S}}
\newcommand{\tV}{\tilde{V}}
\newcommand{\oV}{\overline{V}}

\newcommand{\tM}{\tilde{M}}
\newcommand{\ww}{\left|w\right|}
\newcommand{\piH}{\pi^\ast\negthinspace H}
\newcommand{\piS}{\pi^\ast\negthinspace S}
\newcommand{\piV}{\pi^\ast\negthinspace V}
\newcommand{\Div}{\operatorname{Div}}
\newcommand{\rk}{\operatorname{rk}}
\newcommand{\Pic}{\operatorname{Pic}}

\newcommand{\mult}{\operatorname{mult}}
\newcommand{\Span}{\operatorname{span}}
\newcommand{\cl}{\operatorname{cl}}
\newcommand{\im}{\operatorname{im}}

\newcommand{\supp}{\operatorname{supp}}
\newcommand{\sing}{\operatorname{sing}}
\newcommand{\smooth}{\operatorname{smooth}}
\newcommand{\blow}{\operatorname{Blow}}
\newcommand{\cohom}[2]{H^{#1}\!\left(#2\right)}

\newcommand{\cohomd}[2]{h^{#1}\!\left(#2\right)}

\newcommand{\obundle}[2]{\OO_{#1}\!\left(#2\right)}
\newcommand{\obundletS}[1]{\OO_{\smash{\tS}}\!\left(#1\right)}
\newcommand{\obundles}[1]{\OO_{#1}}
\newcommand{\mydivisor}[3]{\left(#1\piH{#2}E_{#3}\right)\!/2}
\newcommand{\mylinsys }[3]{\left|\mydivisor{#1}{#2}{#3}\right|}
\newcommand{\mybundle }[3]{\obundletS{\mydivisor{#1}{#2}{#3}}}
\newcommand{\mychi    }[3]{\chi\left(\mybundle{#1}{#2}{#3}\right)}
\newcommand{\mycohom  }[4]{\cohom{#1}{\mybundle{#2}{#3}{#4}}}
\newcommand{\mycohomd }[4]{\cohomd{#1}{\mybundle{#2}{#3}{#4}}}
\newtheorem{definition}{Definition}[section]
{\theorembodyfont{\rm}
    \newtheorem{remark}[definition]{Remark}
    \newtheorem{example}[definition]{Example}
}
\newtheorem{lemma}[definition]{Lemma}
\newtheorem{corollary}[definition]{Corollary}

\newtheorem{proposition}[definition]{Proposition}
\newtheorem{theorem}[definition]{Theorem}
\newcommand{\proof}{{\noindent\bf Proof:\ }}
\newcommand{\proofwith}[1]{{\noindent\bf Proof {#1}:\ }}
\newcommand{\proofend}{$\square$\medskip\par}
\newcommand{\romannum}[1]{\setcounter{hlp}{#1}\roman{hlp})}
\begin{document}
\newcounter{hlp}
\begin{center}
    \Large{Minimal even sets of nodes\\}
    \large{S.~Endra\ss\\}
    \large{\today}
\end{center}
\bigskip
\begin{center}
    {\bf Abstract}\\[3ex]
    \begin{minipage}{10cm}\small
        We extend some results on even sets of nodes which have been
        proved for surfaces up to degree $6$ to surfaces up to
        degree $10$. In particular, we give a formula for the
        minimal cardinality of a nonempty even set of nodes.
    \end{minipage}\\[4ex]
\end{center}
\section{Setup}\label{sect:setup}
Let $S\subset\PthreeC$ be a hypersurface of degree $s$ with $\mu$ ordinary
double points (nodes) as its only singularities. Such a surface will be
called a {\em nodal surface} in the sequel. Denote by
$N=\left\{P_1,\ldots,P_\mu\right\}\subset S$ the set of nodes of $S$.
The maximum number of nodes of a nodal surface of degree $d$ is denoted
classically by $\mud$. There is a lot of (old) literature
on nodal surfaces and estimates for $\mud$ (see \cite{endrass}).
For $d=1,2,\ldots,6$
the numbers $\mud$ are $0,1,4,16,31,65$ and for every
$k\in\left\{0,1,\ldots,\mud\right\}$ there exists at least one nodal
surface of degree $d$ with exactly $k$ nodes. In the case of cubic
nodal surfaces ($d=3$), this follows from Cayley's and Schl\"afli's
classification of singular cubic surfaces \cite{cayley},
\cite{schlaefli}. For quartic nodal
surfaces ($d=4$) the fact that $\mu\left(4\right)=16$ is due to
Kummer \cite{kummer}, whereas the construction of arbitrary nodal quartics
goes back to Rohn \cite{rohn}. The first quintic nodal surface ($d=5$) with
31 nodes has been constructed by Togliatti in 1940 \cite{togliatti}. In 1971,
Beauville \cite{beauville} showed that this is in fact the maximal number.
The construction of sextic nodal surfaces ($d=6$) with $1,\ldots,64$
nodes has been given by Catanese and Ceresa
\cite{cataneseceresa}. In 1994, Barth \cite{barth}
constructed a sextic nodal surface with 65 nodes. Shortly afterwards,
Jaffe and Ruberman \cite{jafferuberman}
proved that 65 is the maximal number. Both Beauville
and Jaffe/Ruberman use the code of a nodal surface in their proofs.
This code is a $\FF$ vector space which carries the information
of the low degree contact surfaces of the nodal surface. If a nodal
surface has ``nearly'' $\mu\left(d\right)$ nodes, its code often becomes
accessible.

Let $v\in\NN$  and denote
$\delta\left(v\right)=2\left(v/2-\left\lfloor v/2\right\rfloor\right)$.
This number is 0 if v is even and 1 if v is odd.
We want to study surfaces $V\subset\Pthree$ of degree $v$
with $S.V=2D$ for a (not necessarily smooth or reduced) curve $D$.
In other words, surfaces $V$ which have contact to $S$ along a curve. Let
$\pi\colon\tPthree\rightarrow\Pthree$ be the embedded resolution of
all nodes of $S$. Given such a surface $V$, the proper transforms
of $S$ and $V$ are calculated as
\begin{equation*}
    \tS=\piS-2\sum_{i=1}^\mu E_i
    \quad\text{and}\quad
    \tV=\piV-\sum_{i=1}^\mu \nu_i E_i,
\end{equation*}
where $E_i=\pi^{-1}\left(P_i\right)$ is the exceptional divisor
corresponding to $P_i$ and $\nu_i=\mult\left(V,P_i\right)$ for every
node $P_i\in N$. On the smooth surface $\tS$ we have
$\tV\sim_{lin}2\tD+\sum_{i=1}^\mu\theta_iE_i$, where $\tD$ is the
proper transform of $D$ and the $\theta_i$'s are nonnegative integers.
Let $H\in\Div\left(\Pthree\right)$ be a hyperplane section, then
\begin{equation*}
    2\tD\sim_{lin}v\piH-
    \sum_{i=1}^\mu\left(\nu_i+\theta_i\right)E_i,
\end{equation*}
where $\tD .E_i=\nu_i+\theta_i=\mult\left(D,P_i\right)=\eta_i$. This shows
that in $\Pic\left(\smash{\tS}\right)$ the divisor class
$\left[\delta\left(v\right)\piH+\sum_{\text{$\eta_i$ odd}}E_i\right]$
is divisible by $2$.
This is a remarkable fact, since every $E_i$ is on $\tS$ a smooth,
rational curve with self intersection $-2$. In particular
$E_i\not\sim_{lin}E_j$ for $i\neq j$.

For any set of nodes $M\subseteq N$, let $E_M=\sum_{P_i\in M}E_i$ be the sum
of exceptional curves corresponding to the nodes in $M$.
\begin{definition}\label{definition:even}
    A set $M\subseteq N$ of nodes of $M$ is called {\em strictly even},
    if the cocycle
    class $\cl\left[E_M\right]\in\cohom{2}{\smash{\tS},\ZZ}$
    is divisible by $2$. $M$ is called {\em weakly even}, if the cocycle
    class $\cl\left[\piH+E_M\right]\in\cohom{2}{\smash{\tS},\ZZ}$
    is divisible by $2$. $M$ is called {\em even} if $M$ is 
    strictly or weakly even.
\end{definition}
So the set of nodes $M=\left\{P_i\in N\mid
\text{$\mult\left(D,P_i\right)$ is odd}\right\}$ through which $D$
passes with odd multiplicity is strictly even if
$v$ is even and weakly even if $v$ is odd.
\begin{definition}\label{definition:cut}
    Let $M\subseteq N$ be an even set of nodes of $S$. If
    $V\subset\Pthree$ is a surface with $S.V=2D$ and $M$ is the
    set of nodes of $S$ through which $D$ passes with odd
    multiplicity, we say that {\em $M$ is cut out by $V$ via $D$}.
\end{definition}
Conversely, if $M\subseteq N$ is even, consider the linear system
$\mylinsys{v}{-}{M}$ for
$v\in\NN$ even if $M$ is strictly even and odd if $M$ is
weakly even.
For $v\gg 0$ this linear
system is nonempty by R.R.~and Serre duality. Then for every $v$ such that
$\mylinsys{v}{-}{M}\neq\emptyset$ and for
every divisor $\oD\in\mylinsys{v}{-}{M}$
we can find a surface $V\subset\PP_3$ of degree $v$ which cuts out $M$.
The construction is as follows: $\oD$ is effective, so it
admits a decomposition $\oD=\tD+\sum_{i=1}^\mu \tau_iE_i$ such
that $\tD$ is effective and contains no exceptional component and all
the numbers $\tau_i$ are nonnegative. In particular we have
for all $j\in\left\{1,\ldots,\mu\right\}$ that
\begin{equation*}
    \tD .E_j=\left(\frac{1}{2}\left(v\piH-E_M\right)
    -\sum_{i=1}^\mu\tau_i E_i\right).E_j=
    \left\{\begin{array}{l@{\quad}l}
        2\tau_j & \text{is even if $P_j\not\in M$,}\\
        2\tau_j +1 & \text{is odd if $P_j\in M$.}
    \end{array}\right.
\end{equation*}
But $2\oD\in\left|v\piH-E_M\right|$ on the surface $\tS$, so
$2\oD$ is cut out by a surface
$\oV\in\left|v\piH-E_M\right|$ in $\tPP_3$. Let $V=\pi_\ast\left(\oV\right)$
and $D=\pi_\ast\left(\oD\right)=\pi_\ast\left(\smash{\tD}\right)$,
then by construction
$S.V=2D$ and $\mult\left(D,P_i\right)=\tD .E_i$ for all $i$.
So $M$ is exactly the set of nodes of $S$ through which $D$ passes
with odd multiplicity.
This shows that $M$ is cut out by $V$ via $D$.
Furthermore we see that only nodal surfaces of even degree do admit weakly
even sets of nodes.

If the surface $V$ cuts out an even set of nodes $M$ on $S$ via $D$, then
in general $D$ is not unique with respect to $M$. The set of these contact
curves is parameterized by the linear system
$L_M=\mylinsys{v}{-}{M}$ which is a projective space
of dimension $\mycohomd{0}{v}{-}{M}-1$.
In particular, if $\mycohomd{0}{v}{-}{M}\geq 1$
then there exists a surface of degree $v$ which cuts out $M$. It is funny
to compute these dimensions, though often not possible.

The canonical divisor of $\tS$ is
$\smash{K_{\tS}}\sim_{lin}\left(s-4\right)\piH$.
Define $\binom{n}{k}=0$ for $n<k$, then Riemann Roch for the
bundle $\mybundle{v}{-}{w}$ reads as
\begin{equation*}
    \mychi{v}{-}{w}=\frac{sv}{8}\left(v-2s+8\right)
    +\binom{s-1}{3}+1-\frac{\ww}{4}.
\end{equation*}
The symmetric difference of two strictly even sets of nodes is
strictly even again,
so the set $C_S=\left\{M\subseteq N\mid\text{$M$ is strictly even}\right\}$
carries the natural structure of a $\FF$ vector space sitting inside
$\FF^\mu$. Hence $C_S$ is a binary linear code, which is called
{\em the code of $S$}. The symmetric difference of two weakly even sets of
nodes is strictly even and the symmetric difference of a strictly even set
and a weakly even
set is weakly even. Thus the set
$\oC_S=\left\{M\subseteq N\mid\text{$M$ is even}\right\}$ is a
binary code of dimension $\dim_{\FF}\left(C_S\right)\leq
\dim_{\FF}\left(\oC_S\right)\leq
\dim_{\FF}\left(C_S\right)+1$ sitting also inside $\FF^\mu$.

The elements of $\oC_S$ are called {\em words},
and for every word $w\in \oC_S$ its weight $\left|w\right|$ is its
number of nodes. Let $e_1,\ldots,e_\mu,h$ be the canonical basis of
$\FF^\mu\oplus\FF$ and consider
\begin{equation*}
    \begin{array}{ccccc}
        \FF^\mu & \overset{j}{\longrightarrow} &
        \FF^\mu\oplus\FF & \overset{\lambda}{\longrightarrow} &
            \cohom{2}{\smash{\tS},\FF} \\
              &         & e_i     & \longmapsto     &
            \cl\left[E_i\right] \bmod 2 \\
              &         & h       & \longmapsto     &
            \cl\left[\piH\right]\bmod 2
    \end{array}
\end{equation*}
The projection of $\ker\left(\lambda\right)$ onto the first factor
is nothing but $\oC_S$, and $\ker\left(\lambda\circ j\right) = C_S$.
 If $s$ is even (resp.~odd),
then $\im\left(\lambda\right)$
(resp.~$\im\left(\lambda\circ j\right)$)
is a total isotropic subspace of $\cohom{2}{\smash{\tS},\FF}$ with respect
to the intersection product.
This shows \cite{beauville} that
\begin{align*}
    \dim_{\FF}\left(C_S\right) &
         \geq\mu -\frac{1}{2}b_2\left(\smash{\tS}\right),\\
    \dim_{\FF}\left(\oC_S\right) &
         \geq\mu +1  -\frac{1}{2}b_2\left(\smash{\tS}\right) 
         \quad\text{($s$ even)}.
\end{align*}
The weight of every word $w\in C_S$ is divisible by 4. If
$s=\deg S$ is even, then the weight of every word is divisible by
8 \cite{catanese}.
\subsection{Coding theory}
We recall some definitions and facts from coding theory \cite{lint},
\cite{wall}. Let $C\subseteq\FF^n$ be a linear code and let $e_1,\ldots,e_n$
be the canonical basis of $\FF^n$. $C$ is called {\em even}
if $2\mid\left|w\right|$ for every $w\in C$ and 
{\em doubly even} if   $4\mid\left|w\right|$ for every $w\in C$.
The dual code of $C$ is defined as
\begin{equation*}
    C^\perp = \left\{v\in\FF^n\mid\left< v,w\right>_{\FF}=0\ \forall w\in C
    \right\}.
\end{equation*}
If $C$ is doubly even, then $C\subseteq C^\perp$. Since
$n=\dim_{\FF}\left(C\right)+\dim_{\FF}\left(C^\perp\right)$
we also get $2\dim\left(C\right)\leq n$ with equality iff $C$ is
self dual. For $w\in C$ the support of $w$ is the linear subspace
of $\FF^n$ which is spanned by the ones of $w$, i.e.~
\begin{equation*}
    \supp\left(w\right) = \Span_{\FF}\left\{e_i\mid
        \left<e_i,w\right>=1\right\}.
\end{equation*}
The image of the projection $p_w\colon C\rightarrow\supp\left(w\right)$
is called projection of $C$ onto the support of $w$ and denoted by
$C_w$. Assume that $2d\mid\left|v\right|$ for all $v\in C$ for some
$d\in\NN$. Since $\left|v+w\right|+2\left|v\cap w\right| =
\left|v\right|+\left|w\right|$ and $p_w\left(v\right)=v\cap w$ we see that
$d\mid v'$ for all $v'\in C_w$.
Now the code $C_S$ of the nodal surface $S$ is always doubly even.
If $s=\deg\left(S\right)$ is even, then $\left(C_S\right)_w$ is doubly even
for all $w\in C_S$.

A $\left[n,k,d\right]$-code is a $k$-dimensional linear code
$C\subseteq \FF^n$ with $\left|w\right|\geq d$
for all $w\in C\setminus\left\{0\right\}$.
Many methods have been found to give bounds on $k$
for fixed $n$ and $d$. One of the simplest to apply is the
\begin{theorem}\label{theorem:griesmer}
    (Griesmer bound) For a $\left[n,k,d\right]$ code always
    $n\geq \sum_{i=0}^{k-1} \left\lceil d/2^i\right\rceil$.
\end{theorem}
\subsection{Examples}
The following examples exhibit the trivial and some of the the well known
cases of even sets of nodes \cite{beauville}.
\begin{example}\label{example:cone}
    Let $S$ be a quadratic cone and let $P_1$ be its node. Every line
    $L\subset S$ runs through $P_1$ and there exists exactly one plane
    $H$ with $S.H=2L$. So $\oC_S$ is spanned by $w=\left\{P_1\right\}$
    and $\mycohomd{0}{}{-}{w}=2$.
\end{example}
\begin{example}\label{example:cubic}
    Let $S$ be a cubic nodal surface, then $C_S$ can only be
    non trivial if $S$ has exactly $\mu\left(3\right)=4$ nodes
    $P_1,\ldots,P_4$. But $b_2\left(\smash{\tS}\right)=7$, so
    $\dim_{\FF}\left(C_S\right)\geq 1$. It follows that
    $\dim_{\FF}\left(C_S\right)=1$ and $C_S$ is spanned by
    $w=\left\{P_1,\ldots,P_4\right\}$. But $w$ is
    cut out by a quadric: Riemann-Roch on $\tS$ gives
    $\mychi{2}{-}{w}=3$. From Serre duality we get
    $\mycohomd{2}{2}{-}{w}=\mycohomd{0}{-4}{+}{w}=0$. One easily checks that
    $\mybundle{4}{-}{w}$ is ample, so by Kodaira vanishing also
    $\mycohomd{1}{2}{-}{w}=\mycohomd{1}{-4}{+}{w}=0$.
    This implies that $\mycohomd{0}{2}{-}{w}=3$, so there exists a two
    parameter family of quadric surfaces which cut out $w$.
\end{example}
\begin{example}\label{example:quartic}
    A quartic nodal surface $S$ with $\mu\left(4\right)=16$ nodes is a 
    Kummer surface. Since $b_2\left(\smash{\tS}\right)=22$, we have
    $\dim_{\FF}\left(C_S\right)\geq 5$. On the other hand all nonzero
    words of $C_S$ must have weight $8$ or $16$. So $C_S$ is a
    $\left[16,k,8\right]$ code for some $k\geq 5$. The Griesmer
    bound implies $k\leq 5$, so 
    $C_S$ is a $\left[16,5,8\right]$ code. Every such code has exactly
    one word of weight 16 and 30 words of weight 8. Moreover $C_S$
    is (up to permutation of columns) spanned by the rows of the
    following table.
    \begin{equation*}
        \newcommand{\rr}{\smash{\hspace*{0.06cm}\blacksquare\hspace*{0.06cm}}}
        \newcommand{\nl}{\\\hline}
        \begin{array}{*{16}{|@{}c@{}}|}\hline
      \rr&\rr&\rr&\rr&\rr&\rr&\rr&\rr&   &   &   &   &   &   &   &   \nl
      \rr&\rr&\rr&\rr&   &   &   &   &\rr&\rr&\rr&\rr&   &   &   &   \nl
      \rr&\rr&   &   &\rr&\rr&   &   &\rr&\rr&   &   &\rr&\rr&   &   \nl
      \rr&   &\rr&   &\rr&   &\rr&   &\rr&   &\rr&   &\rr&   &\rr&   \nl
      \rr&\rr&\rr&\rr&\rr&\rr&\rr&\rr&\rr&\rr&\rr&\rr&\rr&\rr&\rr&\rr\nl
        \end{array}
    \end{equation*}
\end{example}
\begin{example}\label{example:quintic}
    A quintic nodal surface $S$ with $\mu\left(5\right)=31$ nodes is
    called Togliatti surface. One computes $b_2\left(\smash{\tS}\right)=53$,
    so again $\dim_{\FF}\left(C_S\right)\geq 5$. By \cite{beauville}, all
    even sets of nodes on $S$ have weight $16$ or $20$. So $C_S$
    is a $\left[31,k,16\right]$ code for some $k\geq 5$. The Griesmer
    bound gives $31\geq 16+8+4+2+1+\left(k-5\right)$, so $k\leq 5$.
    This shows that $C_S$ is a $\left[31,5,16\right]$ code. Every such code
    has exactly 31 words of weight 16 and no word of weight 20.
    Moreover, $C_S$ is (up to a permutation of columns) 
    spanned by the rows of the following table.
    \begin{equation*}
        \newcommand{\rr}{\hspace*{0.03cm}\blacksquare\hspace*{0.03cm}}
        \begin{array}{*{31}{|@{}c@{}}|}\hline
      \rr&\rr&\rr&\rr&\rr&\rr&\rr&\rr&\rr&\rr&\rr&\rr&\rr&\rr&\rr&\rr&
         &   &   &   &   &   &   &   &   &   &   &   &   &   &   \\\hline
      \rr&\rr&\rr&\rr&\rr&\rr&\rr&\rr&   &   &   &   &   &   &   &   &
      \rr&\rr&\rr&\rr&\rr&\rr&\rr&\rr&   &   &   &   &   &   &   \\\hline
      \rr&\rr&\rr&\rr&   &   &   &   &\rr&\rr&\rr&\rr&   &   &   &   &
      \rr&\rr&\rr&\rr&   &   &   &   &\rr&\rr&\rr&\rr&   &   &   \\\hline
      \rr&\rr&   &   &\rr&\rr&   &   &\rr&\rr&   &   &\rr&\rr&   &   &
      \rr&\rr&   &   &\rr&\rr&   &   &\rr&\rr&   &   &\rr&\rr&   \\\hline
      \rr&   &\rr&   &\rr&   &\rr&   &\rr&   &\rr&   &\rr&   &\rr&   &
      \rr&   &\rr&   &\rr&   &\rr&   &\rr&   &\rr&   &\rr&   &\rr\\\hline
        \end{array}
    \end{equation*}
\end{example}
\begin{example}\label{example:sextic}
    Let $S$ be a nodal sextic surface with $\mu\left(6\right)=65$ nodes.
    Every nonzero word $w\in C_S$ must have weight
    $24$, $32$, $40$ or $56$ \cite{jafferuberman}.
    We have $b_2\left(\smash{\tS}\right)=106$,
    so $\dim_{\FF}\left(C_S\right)\geq 12$.

    If $C_S$ contains no word of weight $56$, then 
    $\dim_{\FF}\left(C_S\right)=12$ \cite{jafferuberman},
    \cite{wall}. A short argument runs as follows:
    By the Griesmer bound $C_S$ contains a word $w$ of weight $24$.
    Clearly $p_w\colon C_S\rightarrow\left(C_S\right)_w$ has trivial
    kernel, so $\left(C_S\right)_w$ is a doubly even
    $\left[24,\dim_{\FF}\left(C_S\right),4\right]$ code.
    Hence $\dim_{\FF}\left(C_S\right)\leq 12$.

    It is not clear if $C_S$ is unique up to permutation. It is also 
    not known if any nodal sextic surface can have
    even sets of $56$ or $64$ nodes.
\end{example}
\subsection{The theorem}
For nodal surfaces of degree $6$, Jaffe and Ruberman proved that the smallest
possible nonzero strictly even sets of nodes are the ones cut out by
quadrics. This seems to be true for nodal surfaces of arbitrary
degree, though we only can prove a few cases. For weakly even sets
of nodes, the corresponding statement is proved.
\begin{definition}\label{definition:eminmax}
    For $s\in\NN$ the \emph{minimal cardinality of an even set of nodes
    on a nodal surface of degree $s$} is defined as
    \begin{align*}
        e_{min}\left(s\right) &= \min\left\{
            \left|w\right|\mid\text{$w\in C_S$, $S$ nodal of degree $s$}
        \right\}, \\
        \ose_{min}\left(s\right) &= \min\left\{
            \left|w\right|\mid\text{$w\in \oC_S$, $S$ nodal of degree $s$}
        \right\}.
    \end{align*}
\end{definition}
Our main result is the following
\begin{theorem}\label{theorem:main}
    \begin{itemize}
        \item[\romannum{1}] (Strictly even sets of nodes)
            Let $s\in\left\{3,4,5,6,7,8,10\right\}$. Then
            \begin{equation*}
                e_{min}\left(s\right) =\left\{\begin{array}{cl}
                    s\left(s-2\right) & \text{if $s$ is even,}\\
                    \left(s-1\right)^2 & \text{if $s$ is odd.}
                \end{array}\right.
            \end{equation*}
            Moreover $\left|w\right|=e_{min}\left(s\right)$ if
            and only if
            $w$ is cut out by a quadric surface.
        \item[\romannum{2}] (Weakly even sets of nodes)
            Let $s\in\left\{2,4,6,8\right\}$. Then
            \begin{equation*}
                \ose_{min}\left(s\right)=\frac{s\left(s-1\right)}{2}.
            \end{equation*}
            Moreover $\left|w\right|=\ose_{min}\left(s\right)$ if
            and only if $w$ is cut out by a plane.
    \end{itemize}
\end{theorem}
A close examination of the proof of theorem \ref{theorem:main} exhibits
that certain weights strictly greater than $e_{min}\left(s\right)$ and
$\ose_{min}\left(s\right)$ cannot appear.
\begin{corollary}\label{corollary:main}
    For any nodal surface $S$ of degree $s$, there exist no 
    even sets of nodes with the following weights.
    \begin{equation*}
        \begin{array}{|c||c|c|c|c|}\hline
            s & 6 & 7 & 8 & 10 \\\hline\hline
            \text{weakly even}   & 19,23 & & 32,36,\ldots,56 & \\\hline
            \text{strictly even} &       & 40 & 56 & 88,96,104,112 \\\hline
        \end{array}
    \end{equation*}
\end{corollary}
If $w\in\oC_S$ is cut out by a smooth cubic surface, then
$\ww=3s\left(s-3\right)/2$ \cite{catanese}. The corollary states
that all weights in the open interval
$\left] \ose_{min}\left(s\right),3s\left(s-3\right)/2\right[$ do not
appear for weakly even set of nodes.
In the case of strictly even sets of nodes, the gap is the interval
$\left] e_{min}\left(s\right),2s\left(s-4\right)\right[$. Note that
if $w\in C_S$ is cut out by a smooth quartic surface, then
$\ww=2s\left(s-4\right)$.
\begin{remark}
    It follows from example \ref{example:cone} and example
    \ref{example:cubic} that the theorem is true for $s=2,3$.
\end{remark}
\subsubsection*{Acknowledgments}
    I would like to thank D.~van Straten for valuable discussions.
\section{The formula of Gallarati}
The contact of hypersurfaces in $\PP_r$ along a $r-2$ dimensional variety
has been (to our knowledge) studied first by D.~Gallarati
\cite{gallarati}.
He stated the following
\begin{theorem}\label{theorem:gallarati}
    Let $F_m$, $G_n\subset\PP_r$ be hypersurfaces of degree $m$ and $n$
    with $F_m.G_n=qC$ for some $r-2$ dimensional variety $C$. Assume
    that $F_m$ and $G_n$ have at most double points on $C$. If the
    singular locus of $F_m$ on $C$ (resp.~$G_n$ on $C$) is a $r-3$
    dimensional variety of degree $t$ (resp.~$s$), then
    \begin{equation*}
        q\left(t-s\right) = mn\left(m-n\right).
    \end{equation*}
\end{theorem}
If one allows the surfaces $F_m$ and $G_n$ to have points of
higher multiplicity on $C$, then simple examples show that this number is
dependent on the local geometry. But the philosophy of Gallarati's
theorem is that in the situation of contact of hypersurfaces
the hypersurface of higher degree must have more or harder
singularities on the contact variety than the hypersurface of
lower degree.

We will prove a variant of the above theorem which gives a lower bound
for the size of an even set of nodes. If $S\subset\PP_3$ is a nodal surface
recall that for every even set of nodes $w$ on $S$ there exists a surface
$V\subset\PP_3$ such that $S.V=2D$ and $w$ is just the set of nodes
of $S$ through which $D$ passes with odd multiplicity. We estimate the
number of nodes through which $D$ passes with multiplicity one.

For a slightly more general setup, let $M$ be a smooth projective threefold
and let $S\subset M$ be a nodal surface. Assume that a surface
$V\subset M$ intersects $S$ as $S.V=rD+D'$, $r\geq 2$, for an
irreducible curve $D$ which is not contained in the support of $D'$.
\begin{definition}\label{definition:smooth}
    A node $P$ of $S$ is called $D$-smooth if $P\in D$ and $P$ is a
    smooth point of $V$.
\end{definition}
This definition is justified by the following
\begin{lemma}\label{lemma:smooth}
    Let $P$ be a node of $S$. If $P$ is $D$-smooth, then
    $P$ is a smooth point of $D$. Moreover $r=2$  and
    $P\not\in\supp\left(D'\right)$.
\end{lemma}
\proof There exists a neighborhood $U$ of $P$ in $M$ which is
biholomorphic to some open neighborhood of the origin $\bzero\in\CC^3$,
so it suffices to prove the lemma for two affine hypersurfaces
$S$, $V\subset\CC^3$. We study the intersection with a general
plane through $\bzero$.

Let $L\cong\PP_2$ be the set of all planes $H\subset\CC^3$ through
$\bzero$ and let $T=T_{\bzero}V\in L$ be the tangent plane to $V$
in $\bzero$. Then for all $H\in L\setminus\left\{T\right\}$, the
curve $C_H=V.H$ is smooth in $\bzero$. The set of all planes
$H\in L$ which have contact to the tangent cone $C_{\bzero}S$ of
$S$ in $\bzero$ is parameterized by a smooth conic $Q\subset L$.
For all $H\in L\setminus Q$, the curve $F_H=S.H$ has an ordinary
double point in $\bzero$. While varying $H$ in
$L\setminus\left(Q\cup\left\{T\right\}\right)$, the tangent lines
$T_{\bzero}C_H$ sweep out $T_{\bzero}V$, while the tangent lines to
both branches of $F_H$ in $\bzero$ sweep out $C_{\bzero}S$.
So there exists a plane $\tH\in
L\setminus\left(Q\cup\left\{T\right\}\right)$ such that $T_{\bzero}C_{\tH}$
is not contained in $C_{\bzero}S$. Therefore
$C_{\tH}$ and $F_{\tH}$ meet transversal in $\bzero$, hence on $\tH$
we have local intersection multiplicity
$\left(F_{\tH}.C_{\tH}\right)_{\bzero}=2$. Then of course
\begin{align*}
    2 &= \left(F_{\tH}.C_{\tH}\right)_{\bzero} 
       = \left(\left.S\right|_{\tH}.\left.V\right|_{\tH}\right)_{\bzero}
       = \left(S.V.\smash{\tH}\right)_{\bzero} \\
      &= \left(\left(rD+D'\right).H\right)_{\bzero}
       = r\left(D.\smash{\tH}\right)_{\bzero}+
          \left(D'.\smash{\tH}\right)_{\bzero}.
\end{align*}
Now $\bzero\in D$ implies $\left(D.\smash{\tH}\right)_{\bzero}\geq 1$. Since
$r\geq 2$ we get $\left(D.\smash{\tH}\right)_{\bzero}=1$, $r=2$ and
$\left(D'.\smash{\tH}\right)_{\bzero}=0$. This proves the lemma.\proofend
Now we give the lower bound for the number of $D$-smooth nodes of $S$.
\begin{proposition}\label{proposition:irreducible}
    Assume that $D\not\subseteq\sing\left(V\right)$ and let $\beta$
    be the number
    of singular points of $V$ on $D$ which are smooth points of $S$.
    Then $S$ has at least $D.\left(S-V\right)+\beta$ nodes which are
    $D$-smooth.
\end{proposition}
\proof To prove the theorem we would like to have everything smooth.
There exists a sequence of blowups (embedded resolution of the
singular locus of $S$, $V$ and $D$)
\begin{equation*}
    \tM =M_n\overset{\pi_n    }{\longrightarrow}
    M_{n-1} \overset{\pi_{n-1}}{\longrightarrow}
    \ldots  \overset{\pi_2    }{\longrightarrow}
    M_1     \overset{\pi_1    }{\longrightarrow}
    M_0=M.
\end{equation*}
Let $S_i$, $V_i$ and $D_i$ denote the proper transforms of $S$, $V$ and $D$
with respect to $\pi_i\circ\pi_{i-1}\circ\ldots\circ\pi_1$. We can
define divisors $D_i'$ by
$S_i.V_i=rD_i+D_i'$ with $D_i\not\subset\supp\left(D_i'\right)$,
$1\leq i\leq n$. Moreover we can arrange the maps $\pi_i$ in such
a way that the following conditions hold.
\begin{itemize}
    \item[\romannum{1}] $\pi_1\colon M_1\rightarrow M$
        is the blowup of $M$ in all
        points which are singular for both $S$ and $V$.
    \item[\romannum{2}]
        $\pi=\pi_{n-1}\circ\ldots\circ\pi_2\colon M_{n-1}\rightarrow M_1$
        is the embedded resolution of the singular locus of $V_1$,
        i.e.~every map $\pi_{i+1}$ is a blowup of $M_i$ centered
        in a smooth variety $Z_i\subset M_i$ such that $Z_i$
        is either a point or a smooth curve, $1\leq i\leq n-2$.
    \item[\romannum{3}]
        $\pi_n\colon\tM\rightarrow M_{n-1}$ is the embedded resolution
        of the singularities of $S_{n-1}$ and $D_{n-1}$.
\end{itemize}
Now one has to keep track of the intersection numbers
$D_i.\left(S_i-V_i\right)$ as $i$ increases. We study each of the
three maps separately.

{\bf\romannum{1}} Let $Z=\left\{P_1,\ldots,P_\alpha\right\}=
\sing\left(S\right)\cap\sing\left(V\right)$, then
$M_1=\blow_Z M$. Denote by $E_j=\pi_0^{-1}\left(P_j\right)$ the
exceptional divisor corresponding to $P_j$. The proper transforms are
calculated as
$S_1=\pi_1^\ast S-2\sum_{j=1}^\alpha E_j$ and
$V_1=\pi_1^\ast V-\sum_{j=1}^\alpha m_jE_j$,
where $m_j=\mult\left(V,P_j\right)\geq 2$, $1\leq j\leq\alpha$.
Then the intersection number can be estimated as
\begin{equation*}
    D_1.\left(S_1-V_1\right) = D_1.\pi_1^\ast\left(S-V\right)
        +\sum_{j=1}^\alpha\left(m_j-2\right)D_1.E_j
        \geq D.\left(S-V\right).
\end{equation*}
This is just the information we need, so let us consider the second case.

{\bf\romannum{2}} Every blowup $\pi_{i+1}$ gives rise to an exceptional divisor
$F_{i+1}=\pi_{i+1}^{-1}\left(Z_i\right)$. In the $\left(i+1\right)$-st
step always $S_i$ is smooth in all points of $S_i\cap Z_i$,
whereas $V_i$ is singular in all points of $Z_i$. So the proper transforms
are
\begin{equation*}
    S_{i+1}=\pi_{i+1}^\ast S_i-n_iF_{i+1}
    \quad\text{and}\quad
    V_{i+1}=\pi_{i+1}^\ast V_i-p_iF_{i+1}
\end{equation*}
where $n_i=\mult\left(S_i,Z_i\right)\in\left\{0,1\right\}$ and
$p_i=\mult\left(S_i,Z_i\right)\geq 2$. So this time the intersection
number in question is just
\begin{align*}
    D_{i+1}.\left(S_{i+1}\right. &-\left.V_{i+1}\right) =
    D_{i+1}.\pi_{i+1}^\ast\left(S_i-V_i\right)+
        \left(p_i-n_i\right)D_{i+1}.F_{i+1} \\
    &= D_i.\left(S_i-V_i\right)+
       \left\{\begin{array}{l@{\quad}l}
           \left(p_i-1\right)\mult\left(D_i,Z_i\right) &
           \text{if $Z_i$ is a point,}\\
           p_i\sum_{P\in Z_i\cap D_i}\mult\left(D_i,P\right) &
           \text{if $Z_i$ is a curve,}
        \end{array}\right. \\
    &\geq  D_i.\left(S_i-V_i\right)+ \#\left(Z_i\cap D_i\right).
\end{align*}
But every singularity of $V$ on $D$ outside  the singular locus
of $S$ counts at least once. So by induction
\begin{equation*}
    D_{n-1}.\left(S_{n-1}-V_{n-1}\right)\geq
    D_1.\left(S_1-V_1\right)+\beta\geq
    D.\left(S-V\right)+\beta.
\end{equation*}
where $\beta=\#\left(\left(\sing\left(V\right)\cap D\right)\setminus
\sing\left(S\right)\right)$.

{\bf\romannum{3}} As for the third case we note that $V_{n-1}$ is smooth and
$S_{n-1}$ is nodal with $S_{n-1}.V_{n-1}=rD_{n-1}+D_{n-1}'$. Either
$D_{n-1}\cap\sing\left(S_{n-1}\right)=\emptyset$ or
$D_{n-1}$ contains at least one node of $S$. But then $r=2$ by lemma
\ref{lemma:smooth}
and $D_{n-1}$ is smooth in all nodes of $S_{n-1}$. In both cases,
$S_{n-1}$ and $D_{n-1}$ do not have common singularities.

Let $P_{\alpha+1},\ldots,P_{\alpha+\eta}$ be the nodes of $S_{n-1}$ on
$D_{n-1}$ and let
$P_{\alpha+\eta+1},\ldots,P_{\alpha+\eta+\tau}$ be the remaining nodes
of $S_{n-1}$. Moreover let $E_j=\pi_n^{-1}\left(P_j\right)$,
$\alpha+1\leq j\leq \alpha+\eta+\tau$. But the embedded resolution of the
singularities of $D_{n-1}$ on $V_{n-1}$ is the same as on $S_{n-1}$, so
the proper transforms are
\begin{align*}
    \tS = S_n &= \pi_n^\ast S_{n-1}-E_D-2\sum_{j=\alpha+1}^{\alpha+\eta+\tau}
        E_j, \\
    \tV = V_n &= \pi_n^\ast V_{n-1}-E_D-\sum_{j=\alpha+1}^{\alpha+\eta}E_j-
        \sum_{k=\alpha+\eta+1}^{\alpha+\eta+\tau} q_kE_k,
\end{align*}
where $q_k=\mult\left(V_k,P_k\right)\in\left\{0,1\right\}$ and $E_D$ is a 
sum of exceptional divisors corresponding to the singularities of $D_{n-1}$.
Set $\tD=D_n$ and calculate
\begin{align*}
    \tD.\left(\smash{\tS-\tV}\right) &=
    \tD.\pi_n^\ast\left(S_{n-1}-V_{n-1}\right) -
        \sum_{j=\alpha+1}^{\alpha+\eta}\tD.E_j \\
        &= D_{n-1}.\left(S_{n-1}-V_{n-1}\right)-\eta \\
        &\geq D.\left(S-V\right)-\eta+\beta.
\end{align*}
On the other hand the smooth surfaces $\tS$ and $\tV$ have contact of order
$r-1\geq 1$ along the smooth curve $\tD$. So the tangent bundles
$T_{\tS}$ and $T_{\tV}$ agree along $\tD$. This implies that
the normal bundles $N_{\tD\mid\tS}$ and $N_{\tD\mid\tV}$ coincide, thus
\begin{equation*}
    \left(\smash{{\tD}^2}\right)_{\tS}=
    \deg\left(\smash{N_{\tD\mid\tS}}\right)=
    \deg\left(\smash{N_{\tD\mid\tV}}\right)=
    \left(\smash{{\tD}^2}\right)_{\tV}.
\end{equation*}
Now by adjunction formula $\tD.K_{\tV}=\tD.K_{\tS}$. Using the adjunction
formula again we see that
\begin{equation*}
    0 = \tD.K_{\tS}-\tD.K_{\tV}
      = \tD.\left.\left(\smash{K_{\tM}+\tS}\right)\right|_{\tS}-
         \tD.\left.\left(\smash{K_{\tM}+\tV}\right)\right|_{\tV}
      = \tD.\left(\smash{\tS-\tV}\right).
\end{equation*}
This gives the desired formula $\eta\geq D.\left(S-V\right)+\beta$.
If $V$ is also a nodal surface one can see easily that we have
equality.\proofend
The application to surfaces in $M=\PP_3$ gives the following
\begin{corollary}\label{corollary:irreducible}
    Let a nodal surface $S\subset\PP_3$ of degree $s$ and an irreducible
    surface $V\subset\PP_3$ of degree $v$ intersect as $S.V=2D$ for
    curve $D$ on $S$. Assume that $V$ is not singular along
    a curve contained in $S$ and let $\beta$ be the number of singular
    points of $V$ which are smooth for $S$. If $s>v$, then $D$ is
    reduced. Moreover $V$ cuts out an even set of at least
    $sv\left(s-v\right)/2+\beta$ nodes on $S$ with equality if $V$
    is also nodal.
\end{corollary}
\proof Just run proposition \ref{proposition:irreducible}
on every irreducible component of $D$.\proofend
It is possible to extend proposition  \ref{proposition:irreducible}
to the case when the surface $V$ is not irreducible, but reduced.
The proof however works different.
\begin{proposition}\label{proposition:reduced}
    Let $S\subset\PP_3$ be a nodal surface, $n\in\NN$ and let
    $V_1,\ldots,V_n\subset\PP_3$ be different
    irreducible surfaces of degrees $v_1,\ldots,v_n$ satisfying the
    following conditions:
    \begin{itemize}
        \item[\romannum{1}] $V_i$ is not singular along a curve
            contained in $S$, $1\leq i\leq n$,
        \item[\romannum{2}] $v_i=\deg\left(V_i\right)<s$, $1\leq i\leq n$ and
        \item[\romannum{3}] $S.\left(V_1+\ldots +V_n\right)=2D$ for a (not
            necessarily reduced) divisor $D$ on $S$.
    \end{itemize}
    Then the reduced surface $V=V_1+\ldots +V_n$ of degree
    $v=v_1+\ldots +v_n$ cuts out an even set of nodes $w\in\oC_S$
    of weight $\left|w\right|\geq sv\left(s-v\right)\nx /2$.
\end{proposition}
\proof Since $v_i<s$ and $V_i$ is not singular along a curve contained
in $S$ there exist reduced divisors $D_i$ and $R_i$ on $S$ which
do not have a common component such that $S.V_i=2D_i+R_i$,
$1\leq i\leq n$. But
\begin{equation*}
    S.V=S.\left(V_1+\ldots +V_n\right)=
    2\left(D_1+\ldots +D_n\right)+R_1+\ldots +R_n.
\end{equation*}
This implies that $R_i\subset\bigcup_{j\neq i}V_j$ and thus $R_i$ has
a decomposition $R_i=\sum_{j\neq i}R_{i,j}$ such that
$R_{i,j}\subset V_i\cap V_j$. Now we count the nodes of $S$ through
which $D$ passes with multiplicity 1. Denote $d_i=\deg\left(D_i\right)$,
$r_i=\deg\left(R_i\right)$ and $r_{i,j}=\deg\left(R_{i,j}\right)$.
By corollary \ref{corollary:irreducible},
$V_i$ contains at least $d_i\left(s-v_i\right)$
nodes of $S$ through which $D_i$ passes with multiplicity 1. All these
nodes lie outside $R_i$. We cannot simply add these numbers:
some nodes might be counted more than once. But every node $P\in w$
on $V_i$ which is counted more than once is contained also
in some $V_j$ for a $j\neq i$, hence in $F_{i,j}=V_i.V_j$.
Let  $f_{i,j}=\deg\left(F_{i,j}\right)$. Let $C$ be an irreducible
component of $F_{i,j}$ and let $c=\deg\left(C\right)$. 
We have the following possibilities:
\begin{itemize}
    \item $C\not\subset S$. In this case $C$ contains
        at most $cs/2$ nodes of $S$.
    \item $C\subset S$ is a component of $R_i$. Here $C$ does not contain
        any node that we counted.
    \item $C\subset S$ is a component of $D_i$ and $D_j$. Here $V_i$
        and $V_j$ have contact to $S$ along $C$
        and thus $C$ appears in $F_{i,j}$ with
        multiplicity $\geq 2$. Clearly $C$ contains at most
        $c\left(s-1\right)$ nodes of $S$.
    \item $C\subset S$ is a component of $D_i$, but not of $D_j$.
        Then $V_i$ and $V_j$ meet transversal along $C$, so
        $S.\left(V_i+V_j\right)=3C+other\ curves$. So there exist
        a $k\not\in\left\{i,j\right\}$ such that $C\in F_{i,k}$.
        So $C$ appears with multiplicity $\geq 2$ in
        $\sum_{j\neq i}F_{i,j}$. Again $C$ contains at most
        $c\left(s-1\right)$ nodes of $S$.
\end{itemize}
Since every component of $R_{i,j}$ is contained in $F_{i,j}$, this
shows that $F_{i,j}$ contains at most
$\left(f_{i,j}-r_{i,j}\right) s/2$ nodes that we counted. So $V_i$
contains at least
\begin{equation*}
    d_i\left(s-v_i\right)-\sum_{j\neq i}\left(f_{i,j}-r_{i,j}\right)
\end{equation*}
nodes through which $D$ passes with multiplicity 1. This implies
\begin{align*}
    \left|w\right| &\geq
    \sum_{i=1}^n d_i\left(s-v_i\right)
        -\frac{s}{2}\sum_{j\neq i}\left(f_{i,j}-r_{i,j}\right)\\
    &=\frac{1}{2}\left(\sum_{i=1}^n\left(sv_i-r_i\right)\left(s-v_i\right)
        +sr_i -s\sum_{j\neq i}v_iv_j\right)\\
    &=\frac{1}{2}\left(s^2\left(v_1+\ldots v_n\right)
        -s\left(v_1^2+\ldots +v_n^2+2\sum_{i<j}v_iv_j\right)
        +r_1v_1+\ldots +r_nv_n\right)\\
    &\geq \frac{sv}{2}\left(s-v\right)
\end{align*}
This completes the proof.\proofend
\section{Contact surfaces and quadratic systems}
In this section we apply the previous results to our initial
situation. So let again $S\subset\PP_3$ be a nodal surface of degree $s$ and
$V\subset\PP_3$ a reduced surface of degree $v$ such that $S.V=2D$ for some
curve $D$. We give a complete analysis of
the situation when $V$ is a plane or a quadric.
Using the notation of the first paragraph, $V$ cuts out an even
set of nodes $w\in\oC_S$. Recall that the linear system
$L_w=\mylinsys{v}{-}{w}$ parameterizes all contact curves of the form
$D'=(1/2)S.V'$ where $V'$ is a surface of degree $v$ which cuts out $w$.
In some cases, $V$ will be the unique surface
of degree $v$ which cuts out $w$.

Now $2D\in\PP\left(\cohom{0}{\obundletS{vH}}\right)$ is the restriction of
$V\in\PP\left(\cohom{0}{\obundle{\PP_3}{vH}}\right)$ to $S$. Consider the
exact sequence
\begin{equation*}
    0\longrightarrow
    \obundle{\PP_3}{\left(v-s\right)H}\longrightarrow
    \obundle{\PP_3}{vH}\longrightarrow
    \obundle{S}{vH}\longrightarrow
    0.
\end{equation*}
Since $\cohom{i}{\obundle{\PP_3}{\left(v-s\right)H}}=0$ for $s>v$,
$i=0,1$, the induced map
$\cohom{0}{\obundle{\PP_3}{vH}}\rightarrow
\cohom{0}{\obundle{S}{vH}}$ is an isomorphism. So if $v<s$, then $V$ is the
unique surface of degree $v$ cutting out $w$ via $D$ and
$L_w=\mylinsys{v}{-}{w}$ in fact parameterizes the
space of all surfaces of degree $v$ which cut out $w$. This space is
not a linear system, but the quadratic system
\begin{equation*}
    Q_w=\left\{V'\mid S.V'=2D'\text{\ with\ }D'\in L_w\right\}.
\end{equation*}
It is constructed as follows:
If $\mycohomd{0}{v}{-}{w}=n+1\geq 2$ we can
find $n+1$ linearly independent sections
$s_0,\ldots,s_n\in\mycohom{0}{v}{-}{w}$.
Clearly all products $s_is_j\in\cohom{0}{\obundletS{v\piH-E_w}}$
for $0\leq i\leq j\leq n$.
So there exist sections $g_{i,j}\in\cohom{0}{\obundle{\tPP_3}{v\piH-E_w}}$
over $\tPP_3$ which restrict to $s_is_j$ under the
identification $\obundle{\tPP_3}{v\piH-E_w}\otimes \obundles{\tS}\cong
\obundletS{v\piH-E_w}$.
% But $L_w=
%\left\{\lambda_0s_0+\ldots+\lambda_ns_n=0\mid
%\left(\lambda_0:\ldots :\lambda_n\right)\in\PP_n\right\}$. Moreover,
Outside the exceptional locus we can view the $g_{i,j}$ as sections
of $\obundle{\PP_3}{vH}$. Since $w$ has codimension $\geq 2$ in
$\PP_3$, these sections extend also to $w$. This implies that
\begin{align*}
    Q_w &= Q\left(g_{i,j}\mid 0\leq i\leq j\leq n\right) \\
        &= \left\{\sum_{i=0}^n\lambda_i^2g_{i,i}+
    2\sum_{0\leq i<j\leq n}\lambda_i\lambda_j g_{i,j}=0\mid
    \left(\lambda_0:\ldots:\lambda_n\right)\in\PP_n\right\}.
\end{align*}
Therefore the quadratic system $Q_w$ is the image of an embedding of
Veronese type of $\PP_n$ into the space
$\PP_{\binom{v+3}{3}}$ parameterizing all surfaces of degree $v$.
In general, $Q_w$ will not contain any linear subspace.

The quadratic system $Q_w$ admits a decomposition $Q_w=B_w+F_w$ where
$B_w$ is a reduced surface of degree $b\leq v$ and the base locus of
$F_w$ (if any) consists only of curves and points. If $F_w$ has no
basepoints then $B_w$ cuts out $w$ and so
$\mycohomd{0}{b}{-}{w}=1$.
\begin{definition}\label{definition:stable}
    An even set of nodes $w\in\oC_S$ is called
    \begin{equation*}
        \left.\begin{array}{c}
            \text{semi stable}\\\text{stable}\\\text{unstable}
        \end{array}\right\}
        \quad\text{in degree $v$ if}\quad
        \left\{\begin{array}{c}
            \text{$F_w$ is basepointfree,}\\
            \text{$F_w=\emptyset$,}\\
            \text{$F_w$ has basepoints.}
        \end{array}\right.
    \end{equation*}
\end{definition}
The base locus of $Q_w$ is
$B\left(Q_w\right)=\left\{g_{i,j}=0\mid 0\leq i\leq j\leq n\right\}$.
It is contained in the discriminant locus
$Z\left(Q_w\right)=\left\{g_{i,i}g_{j,j}=g_{i,j}^2
\mid 0\leq i<j\leq n\right\}$. There is also a Bertini type theorem
for quadratic systems.
\begin{lemma}\label{lemma:bertini}
    (Bertini for quadratic systems) The general element of $Q_w$
    is smooth outside $Z\left(Q_w\right)$.
\end{lemma}
\proof The proof runs like the proof of the Bertini theorem in
\cite{griffithsharris}.\proofend
Next we give a different characterization of stability.
\begin{proposition}\label{proposition:stable}
    Let $w\in\oC_S$. The surface $B_w$ is always reduced and
    \begin{itemize}
        \item[\romannum{1}] $w$ is stable in degree $v$ if and only if
            $\mycohomd{0}{v}{-}{w}=1$.
        \item[\romannum{2}] $w$ is semi stable in degree $v$ if and only if
            $F_w$ contains a square. Then $B_w$ cuts out
            $w$ and either every surface in $F_w$ is a square
            or the general surface in $F_w$ is reduced.
        \item[\romannum{3}] $w$ is unstable in degree $v$ if and only if
            $F_w$ contains no square. Then $B_w$ does not
            cut out $w$ and the general surface in $Q_w$ is reduced.
    \end{itemize}
\end{proposition}
\proof \romannum{1} follows from the definition. So let $w$ be not stable
in degree $v$. We use induction on $n$.

$n=1$: By construction $\gcd\left(g_{0,0},g_{0,1},g_{1,1}\right)=g$ is
reduced. Let $\og_{i,j}=g_{i,j}/g$ and let 
$\oQ_w=Q\left(\og_{0,0},\og_{0,1},\og_{1,1}\right)$. Now we have two cases.

a) If $\og_{0,0}\og_{1,1}=\og_{0,1}^2$, then $\og_{0,0}$ and $\og_{1,1}$
must be squares. So $\og_{0,0}=a_0^2$, $\og_{1,1}=a_1^2$ and thus
$\og_{0,1}=a_0a_1$. Hence $B_{w}=\left\{g=0\right\}$ cuts out $w$. So
by  construction the quadratic system 
$F_w=\left\{\smash{\left(\lambda_0a_0+\lambda_1a_1\right)^2}
\mid\left(\lambda_0:\lambda_1\right)\in\PP_1\right\}$ contains only
squares. Then $F_w$ is free and $w$ is semi stable in degree $v$.

b) $Z\left(\oQ_w\right)= \left\{\og_{0,0}\og_{1,1}=\og_{0,1}^2\right\}$
is a surface. If all surfaces of $Q_w$ are not reduced,
then by lemma \ref{lemma:bertini} all surfaces of
$\oQ_w$ contain a component of $Z\left(Q_w\right)$.
So this component is constant for all surfaces in  $\oQ_w$,
which contradicts
$\gcd\left(\og_{0,0},\og_{0,1},\og_{1,1}\right)=1$.  So the general surface
in $Q_w$ is reduced. Now assume $B_w$ cuts out $w$. Then by construction
$F_w$ contains squares, so $F_w$ is free and $w$ is semi stable in
degree $v$. Otherwise $B_w$ does not cut out $w$, so $F_w$ must have
basepoints in $w$. Then $w$ is unstable in degree $v$.

$n-1\Rightarrow n$: Again
$\gcd\left(g_{i,j}\mid 0\leq i\leq j\leq n\right)=g$ is reduced. Consider
the quad\-ra\-tic system
$Q=Q\left(g_{i,j}\mid 0\leq i\leq j\leq n-1\right)$. Either
$\gcd\left(g_{i,j}\mid 0\leq i\leq j\leq n-1\right)$ is reduced and
we're done or it's not reduced. For
$\lambda=\left(\lambda_0:\ldots:\lambda_{n-1}\right)\in\PP_{n-1}$
let
\begin{equation*}
    g_\lambda =\sum_{i=0}^{n-1}\lambda_i^2g_{i,i}
    +2\sum_{0\leq i<j\leq n-1}\lambda_i\lambda_j g_{i,j}
    \quad\text{and}\quad
    h_\lambda=\sum_{i=0}^{n-1}\lambda_ig_{i,n}.
\end{equation*}
Now consider the quadratic system
\begin{equation*}
    R_\lambda =\left\{
    t^2 g_\lambda+2t\lambda_nh\lambda+\lambda_n^2g_{n,n}=0
    \mid\left(t:\lambda_n\right)\in\PP_1\right\}.
\end{equation*}
While varying $\lambda\in\PP_{n-1}$, 
$\gcd\left(g_\lambda,h_\lambda,g_{n,n}\right)$ is constant on an open
dense subset, since it contains only factors of $g_{n,n}$. So for
general $\lambda$ $\gcd\left(g_\lambda,h_\lambda,g_{n,n}\right)=
\gcd\left(g_{i,j}\mid 0\leq i\leq j\leq n\right)=g$ is reduced.
By the first part either $g_\lambda g_{n,n}=h_\lambda^2$ for all
$\lambda$, so $g_\lambda$ and $g_{n,n}$ are always squares modulo
$g$. Then $w$ is semi stable in degree $v$. Or the general surface
in $R_\lambda$ and hence in $Q_w$ is reduced. Again either
$B_w$ cuts out $w$ and we're in the semi stable case or $B_w$ does not cut
out $w$. Then $w$ is unstable in degree $v$.\proofend
\begin{corollary}\label{lemma:unique}
    If $2v<s$ then $w$ is semi stable in degree $v$.
\end{corollary}
\proof On $S$ we have 
$g_{0,0}g_{1,1}-g_{0,1}^2=s_0^2s_1^2-\left(s_0s_1\right)^2=0$. Let
$S=\left\{f=0\right\}$, then either $g_{0,0}g_{1,1}=g_{0,1}^2$ or
$f\mid g_{0,0}g_{1,1}-g_{0,1}^2$. But the second case implies
$2v\geq s$, so we are in the first case. Then we always run
into case a) in the proof of proposition \ref{proposition:stable}.\proofend
\begin{corollary}
    If $w$ is semi stable in degree $v$ and $2\deg\left(F_w\right)<s$,
    then $F_w$ contains only squares.
\end{corollary}
\proof $F_w$ contains a square $W_0=\left\{g_{0,0}=g_0^2=0\right\}$.
Now take any other $W_1=\left\{g_{1,1}=0\right\}\in F_w$ and consider
the quadratic system generated by $g_{0,0}$ and $g_{1,1}$:
$g_{0,0}$, $g_{1,1}$ give rise to sections $s_0^2$, $s_1^2$ over
$\tS$. Then $s_0s_1$ is the restriction of a section $g_{0,1}$ to
$\tS$. The quadratic system in question is just
$Q=Q\left(g_{0,0},g_{0,1},g_{1,1}\right)$. But 
$g_{0,0}g_{1,1}-g_{0,1}^2$ vanishes on $S$. Since
$\deg\left(g_{0,0}g_{1,1}-g_{0,1}^2\right)=2\deg\left(F_w\right)<s$,
we have $g_{0,0}g_{1,1}=g_0^2g_{1,1}=g_{0,1}^2$. This implies that
also $g_{1,1}$ is a square.\proofend
\begin{proposition}\label{proposition:unstable}
    If $w$ is unstable in degree $v$, then there exists a surface
    $W$ of degree $2v-s$ such that $w$ is cut out by a reduced surface
    $V$ of degree $v$ satisfying:
    \begin{itemize}
        \item[\romannum{1}] $V$ is not singular on $S$ outside $W$.
        \item[\romannum{2}] If $V$ is singular along a curve $C\subset S$,
            then $C$ is a curve of triple points of $W$.
    \end{itemize}
\end{proposition}
\proof $w$ is cut out by a reduced surface, so we can assume that
$g_{0,0}$ is square free.
Again $g_{0,0}g_{1,1}-g_{0,1}^2$ vanishes on $S$ and $g_{0,0}$,
$g_{1,1}$ are linearly independent. So there exists a polynomial
$\alpha$ of degree $2v-s$ such that $\alpha f=g_{0,0}g_{1,1}-g_{0,1}^2$.
Let $W=\left\{\alpha =0\right\}$ and let
$V_\lambda=\left\{
\lambda_0^2 g_{0,0}+2\lambda_0\lambda_1 g_{0,1}+\lambda_1^2 g_{1,1}=0
\mid\lambda=\left(\lambda_0:\lambda_1\right)\in\PP_1\right\}$.
For every point $P\in\PP_3$ we can choose affine coordinates
$\left(z_1,z_2,z_3\right)$ on an affine neighborhood $U$ of $P$.
For any function $h$ on $U$, we identify the total derivative
$Dh$ with the gradient $\nabla h$ and $D^2h$ with the Hesse matrix
$H\left(h\right)$. We find that
\begin{align*}
    D\left(\alpha f\right) &=
        \alpha\nabla f+f\nabla\alpha =
        g_{0,0}\nabla g_{1,1}+g_{1,1}\nabla g_{0,0}-
       2g_{0,1}\nabla g_{0,1},\\
    D^2\left(\alpha f\right) &=
        \alpha H\left(f\right)+f H\left(\alpha\right)+
        \nabla\alpha{\nabla f}^t+\nabla f{\nabla\alpha}^t\\
        &= g_{0,0}H\left(g_{1,1}\right)+
           g_{1,1}H\left(g_{0,0}\right)-
          2g_{0,1}H\left(g_{0,1}\right)\\
        &\hspace*{4ex}+
           \nabla g_{0,0}{\nabla g_{1,1}}^t+
           \nabla g_{1,1}{\nabla g_{0,0}}^t-
          2\nabla g_{0,1}{\nabla g_{0,1}}^t.
\end{align*}
Now let $P\in S$. We have to consider two different cases:

a) $P\in\sing\left(S\right)\setminus W$, so $f\left(P\right)=0$,
$\nabla f\left(P\right)=0$ and 
$\rk\left(H\left(f\right)\left(P\right)\right)=3$.
If $P$ is a basepoint of $Q$ then
\begin{align*}
    H\left(\alpha f\right)\left(P\right) &=
        \alpha\left(P\right)H\left(f\right)\left(P\right) \\
    &= \left(\nabla g_{0,0}{\nabla g_{1,1}}^t+
           \nabla g_{1,1}{\nabla g_{0,0}}^t-
          2\nabla g_{0,1}{\nabla g_{0,1}}^t\right)\left(P\right).
\end{align*}
But $P\not\in W$, so $\alpha\left(P\right))\neq 0$ and
$\rk\left(H\left(\alpha f\right)\left(P\right)\right)=3$. This is only
possible if $\nabla g_{0,0}$, $\nabla g_{1,1}$ and $\nabla g_{0,1}$
are linearly independent in $P$. So every surface $V_\lambda$
is smooth in $P$. If $P$ is not a basepoint of $Q$ then the general
surface $V_\lambda$ will not contain $P$.

b) Let $P\in\smooth\left(S\right)\setminus W$. Here
$f\left(P\right)=0$, $\nabla f\left(P\right)\neq 0$ and
$\alpha\left(P\right)\neq 0$. Then
$\nabla\left(\alpha f\right)\left(P\right)=
\alpha\left(P\right)\nabla f\left(P\right)\neq 0$, so $P$ is not a basepoint
of $Q$. Assume now we have chosen $\lambda$
such that $P\in V_\lambda$. After a permutation
of indices we can assume $\lambda_0=1$, so
$V_\lambda=\left\{
g_{0,0}+2\lambda_1 g_{0,1}+\lambda_1^2 g_{1,1}=0\right\}$.
Since $P$ is not a basepoint we have $g_{1,1}\left(P\right)\neq 0$.
Together with $\left(g_{0,0}g_{1,1}-g_{0,1}^2\right)\left(P\right)=0$
we get $\lambda_1=-\left(g_{0,1}/g_{1,1}\right)\left(P\right)$.
Then
\begin{align*}
    \nabla &\left(g_{0,0}+2\lambda_1 g_{0,1}+\lambda_1^2
    g_{1,1}\right)\left(P\right)=\\
    &\hspace*{5ex}=\frac{1}{g_{1,1}\left(P\right)}
    \left(g_{1,1}\nabla g_{0,0}-2 g_{0,1}\nabla g_{0,1}
    +g_{0,0}\nabla g_{1,1}\right)\left(P\right)\\
    &\hspace*{5ex}=\frac{1}{g_{1,1}\left(P\right)}
    \alpha\left(P\right)\nabla f\left(P\right).
\end{align*}
We see that $P$ is a smooth point of $V_\lambda$,
so together with a) we have proved \romannum{1}.

c) Assume that $V_\lambda$ is singular along a curve
$C_\lambda\subset S$ and let $m_\lambda=\mult\left(V_\lambda,C_\lambda\right)$.
Then $C_\lambda$ is a continuous family of curves and
$m=\min\left\{m_\lambda\mid\lambda\in\PP_1\right\}$ is equal
to $m_\lambda$ on an open dense subset of $\PP_1$. Now
\romannum{1} says that 
$C_\lambda\subseteq S\cap W$ for all $\lambda\in\PP_1$.
But $S\cap W$ is itself a curve, so this family is
in fact constant. So let $C=C_{(0:1)}$. 
Now $g_{0,0}g_{1,1}-g_{0,1}^2=\alpha f$ vanishes to the $2m$-th order
along $C$ and $\mult\left(S,C\right)=1$, so $\alpha$ vanishes
to the $\left(2m-1\right)$-st order along $C$.\proofend
\begin{corollary}\label{corollary:unstable}
    Let $w\in\oC_S$.
    \begin{itemize}
        \item[\romannum{1}]
            If $w$ is unstable in degree $s/2$, then $\ww=s^3/8$.
        \item[\romannum{2}]
            If $w$ is unstable in degree $\left(s+1\right)\!/2$
            (resp.~$\left(s+2\right)/2$),
            then $\ww\geq s\left(s-1\right)^2/8$
            (resp.~$s\left(s-2\right)^2/8$).
    \end{itemize}
\end{corollary}
\proof In the first case $W=\emptyset$. So the general surface in
$Q_w$ is not singular on $S$, hence irreducible. Now apply
corollary \ref{corollary:irreducible}.

In the second case
$\deg\left(W\right)\leq 2$, so $W$ has no triple curve.
Now apply proposition \ref{proposition:reduced}.\proofend
Now here comes our analysis what happens if $V$ is a plane or a quadric.
\begin{proposition}\label{proposition:planequadric}
    Let $w\in\oC_S$.
    \begin{itemize}
        \item[\romannum{1}] If $w$ is cut out by a plane $H$, then
            $\ww=s\left(s-1\right)\!/2$.
            Moreover $w$ is stable in degree
            $1$ if $s>2$ and unstable in degree $1$ otherwise.
        \item[\romannum{2}] If $w$ is cut out by a reduced quadric
            $Q$, then
            \begin{equation*}
                \ww = \left\{\begin{array}{c@{\quad}l}
                    s\left(s-2\right)  & \text{if $s$ is even,}\\
                    \left(s-1\right)^2 & \text{if $s$ is odd.}
                \end{array}\right.
            \end{equation*}
            Moreover $w$ is stable in degree $2$ if $s>4$ and unstable in
            degree $2$ otherwise.
    \end{itemize}
\end{proposition}
\proof \romannum{1} $H$ is smooth, so $\ww=s\left(s-1\right)/2$ by
corollary \ref{corollary:irreducible}.
If $s>2$ then $2\deg\left(H\right)=2<s$, so
$w$ is semi stable in degree $1$ by lemma \ref{lemma:unique}. But
then $w$ is stable in degree $1$.
In the case $s=2$ example \ref{example:cone} shows
that $w$ is unstable in degree $1$.

\romannum{2} Assume first $Q$ is nodal. Then
$\ww\in\left\{s\left(s-2\right),s\left(s-2\right)+1\right\}$ by corollary
\ref{corollary:irreducible}.
But $s\left(s-2\right)+1=\left(s-1\right)^2$ and
$4\mid\ww$ imply the above formula for $\ww$.

Now let $Q=H_1+H_2$ where $H_1\neq H_2$ are planes and
set $L=H_1\cap H_2$. If $L\subset S$, then 
$S.H_i=2D_i+L$ and each $D_i$ contains exactly
$\left(s-1\right)^2\!/2$ nodes which are $D_i$-smooth by proposition
\ref{proposition:irreducible}.
Clearly $L$ cannot contain any node of $w$. So $\ww=\left(s-1\right)^2$.
If $L\not\subset S$, then $S.H_i=2D_i$. Every
$D_i$ is reduced and contains exactly $s(s-1)/2$ $D_i$-smooth nodes.
In every point $P\in L\cap S$
both $H_1$ and $H_2$ are tangent to $S$, so $P$ is a
node of $S$. Both $H_1$ and $H_2$ have contact to the tangent
cone $C_PS$ of $S$ at $P$. This implies $L\not\subset C_PS$,
hence $\mult\left(S,L;P\right)=2$. Therefore $L$ contains exactly
$s/2$ such nodes and $\ww=s\left(s-2\right)$.

If $s>4$, then $w$ is stable in degree $2$. Now let $s\leq 4$.
In any case $F_w$ cannot contain a square. So $w$ is unstable
if $\mycohomd{0}{2}{-}{w}>1$. For $s=3$ this follows from
example \ref{example:cubic}. If $s=4$ then we find using Serre duality that
$\mycohomd{2}{2}{-}{w}=\mycohomd{0}{-2}{+}{w}=0$.
Therefore it follows that 
$\mycohomd{0}{2}{-}{w}\geq\mychi{2}{-}{w}=2$.\proofend
\section{The proof of theorem \ref{theorem:main}}
This section is devoted entirely to the proof of theorem \ref{theorem:main}.
Let $S$ and $V$ with $S.V=2D$ as in the first section.
\begin{lemma}\label{lemma:equal}
    \cite{catanese}
    Let $w\in \oC_S$ be an even set of nodes.
    Let $n\in\NN$ be even if $w$ is strictly even and
    odd if $w$ is weakly even. Then for all $i\geq 0$
    \begin{equation*}
        \cohomd{i}{\obundletS{\left(n\piH +E_w\right)\nx /2}} =
        \cohomd{i}{\obundletS{\left(n\piH -E_w\right)\nx /2}}.
    \end{equation*}
\end{lemma}
\proofwith{ of theorem \ref{theorem:main}}
\romannum{1} By proposition \ref{proposition:planequadric}
the even sets $w\in\oC_S$ cut
out by planes satisfy $\ww=s\left(s-1\right)\!/2$. We show that no smaller
even sets can occur. The proof also explains the ``gaps'' of
corollary \ref{corollary:main}. So let $w\in\oC_S\setminus\left\{0\right\}$ be weakly even.

$s=4$: Since $\mychi{}{-}{w}=\left(10-\ww\right)\!/4$ is an integer
we must have
$\ww\in\left\{2,6,10,14\right\}$. By Serre duality and lemma
\ref{lemma:equal}
$\mycohomd{2}{}{-}{w}=\mycohomd{0}{-}{-}{w}=0$. Now let $\ww\leq 6$, then
$\mycohomd{0}{}{-}{w}\geq 1$. So $w$ is cut out by a plane, hence $\ww=6$.

$s=6$: Here 
$\mychi{}{-}{w}=\mychi{3}{-}{w}=\left(35-\ww\right)\!/4$,
so $\ww\in\left\{3,7,11,\ldots\right\}$. Let $\ww<27$,
 so $\ww\leq 23$ and $\mychi{}{-}{w}\geq 3$. Following proposition
\ref{proposition:stable},
we see that $w$ is either stable in degree $1$ or no plane cuts out $w$,
hence $\mycohomd{2}{3}{-}{w}=\mycohomd{0}{}{-}{w}\in\left\{0,1\right\}$.
This implies $\mycohomd{0}{3}{-}{w}\geq 2$. Now corollary
\ref{corollary:unstable} tells us
that $w$ is semi stable in degree $3$. But $w$ cannot be stable in degree
$3$, so $w$ is stable in degree $1$. Hence
$w$ is cut out by a plane. It follows that $\ww=15$ and that there are
no weakly even sets of $19$ and $23$ nodes on a nodal sextic surface.

$s=8$: Now $\mychi{3}{-}{w}=\mychi{5}{-}{w}=21-\ww /4$, so
$\ww\in\left\{4,8,12,\ldots\right\}$. Let $\ww<60$, then
$\ww\leq 56$ and $\mychi{3}{-}{w}\geq 7$.
Assume that $\mycohomd{0}{3}{-}{w}=\mycohomd{2}{5}{-}{w}\leq 1$. It
follows that $\mycohomd{0}{5}{-}{w}\geq 6$,
hence $w$ is unstable in degree $5$.
So $\ww\geq 60$ by corollary \ref{corollary:unstable}, contradiction.
This implies that $w$ is semi stable in degree $3$ and stable in degree $1$.
Again $w$ is cut out by a plane, hence $\ww=28$. Moreover, a nodal octic
surface cannot have weakly even sets of $32, 36, 40,\ldots,56$ nodes.

\romannum{2} This is essentially a copy of the methods of \romannum{1}.
Let $w\in C_S\setminus\left\{0\right\}$ be strictly even.

$s=4$: We have $\ww\in\left\{8,16\right\}$. If $\ww=9$ then
$\mychi{2}{-}{w}=2$. But $\mycohomd{2}{2}{-}{w}=\mycohomd{0}{-2}{-}{w}=0$,
so by Serre duality and lemma \ref{lemma:equal} 
$\mycohomd{0}{2}{-}{w}\geq 2$ and $w$ is cut out by a quadric.

$s=5$: Now $\ww\in\left\{4,8,12,\ldots\right\}$ and
$\mychi{2}{-}{w}=5-\ww/4$. As usual we find that
$\mycohomd{2}{2}{-}{w}=\cohomd{0}{\obundletS{-E_w/2}}=0$. So if
$\ww\leq 16$ then $\mycohomd{0}{2}{-}{w}\geq 1$ and $w$ is cut out
by a quadric. Then $\ww=16$ by proposition \ref{proposition:planequadric}.

$s=6$: Here $\ww\in\left\{8,16,24,\ldots\right\}$ and
$\mychi{2}{-}{w}=8-\ww /4$. This time we find that
$\mycohomd{2}{2}{-}{w}=\mycohomd{0}{2}{-}{w}$. So if $\ww\leq 24$, then
$\mycohomd{0}{2}{-}{w}\geq 1$ and $w$ is cut out by a quadric surface.
But then $\ww=24$.

$s=7$: We modify the proof as follows. One calculates
$\mychi{2}{-}{w}=\mychi{4}{-}{w}=14-\ww/4$. Let $\ww<44$, so
$\ww\leq40$. By proposition \ref{proposition:planequadric}
$\mycohomd{0}{2}{-}{w}\in\left\{0,1\right\}$, so
$\mycohomd{0}{4}{-}{w}\geq 3$. If $w$ is unstable in degree $4$ then
$\ww\geq42$, contradiction. So $w$ is semi stable in degree $4$ and
stable in degree $2$. Now $w$ is cut out by a quadric and thus
$\ww=36$. In particular, there is no even set of $40$ nodes on
a nodal septic surface.

$s=8$: Using the same argument as for $s=7$, we get
$\mychi{4}{-}{w}=20-\ww/4$ and
$\mycohomd{2}{4}{-}{w}=\mycohomd{0}{4}{-}{w}$. Let $\ww<64$, then
$\ww\leq56$ and $\mycohomd{0}{4}{-}{w}\geq 3$. Again by corollary
\ref{corollary:unstable}
$w$ cannot be unstable in degree $4$. Hence $w$ is semi stable in degree
$4$ and stable in degree $2$. In particular $\ww=48$ and there is
no even set of $56$ nodes on a nodal octic surface.

$s=10$: Finally we calculate $\mychi{6}{-}{w}=40-\ww/4$ and as before
$\mycohomd{2}{6}{-}{w}=\mycohomd{0}{6}{-}{w}$. Let
$\ww<120$, so $\ww\leq 112$ and $\mycohomd{0}{6}{-}{w}\geq 6$. As before
$w$ is semi stable in degree $6$. If $w$ was stable in degree $4$ then
$\mycohomd{0}{6}{-}{w}=4$, contradiction. So $w$ is semi stable in degree
$w$ and stable in degree $2$. Again $\ww=80$, hence there are no
strictly even sets of $88$, $96$, $104$ and $112$ nodes
on a nodal surface of degree $10$.\proofend
\section{Examples revisited}
We want to go a little more into the example of quartics. Many of the
facts stated in this example can be found in \cite{gallarati}.
\begin{example}
    Let $S$ be a nodal quartic surface and let $w\in C_S$ with
    $\ww=8$. We have seen that $\mycohomd{0}{2}{-}{w}\geq 2$
    and that $w$ is unstable in degree $2$. Let $Q_w$ be the
    quadratic system of quadrics which cut out $w$. The base locus
    of $Q_w$ is contained in the surface $W$ of proposition
    \ref{proposition:unstable}. But here $W=\emptyset$, so the
    only basepoints of $Q_w$ are the nodes of $S$. It follows from
    lemma \ref{lemma:smooth} and Bertini that the general element
    in $L_w=\mylinsys{2}{-}{w}$ is a smooth elliptic curve.
    In fact $Q_w$ defines an elliptic fibration of
    $\tS$. If $\mycohomd{0}{2}{-}{w}>2$, then we will find two
    smooth elliptic curves on $\tS$ intersecting in at least one point,
    contradiction. This shows that $\mycohomd{0}{2}{-}{w}=2$.
    Taking into account that $\mydivisor{4}{-}{w}$ is nef and
    big on $\tS$, we can calculate the numbers of the next table.
    If $\ww=16$ one computes
    $\mycohomd{i}{4}{-}{w}$ in the same fashion.
    \begin{equation*}
        \begin{array}{|c||c|c|c|}\hline
            \ww=8               & h^0 & h^1 & h^2 \\\hline\hline
            \mydivisor{2}{-}{w} &  2  &  0  &  0  \\\hline
            \mydivisor{4}{-}{w} &  8  &  0  &  0  \\\hline
        \end{array}
        \quad\quad   
        \begin{array}{|c||c|c|c|}\hline
            \ww=16              & h^0 & h^1 & h^2 \\\hline\hline
            \mydivisor{2}{-}{w} &  0  &  0  &  0  \\\hline
            \mydivisor{4}{-}{w} &  6  &  0  &  0  \\\hline
        \end{array}
    \end{equation*}
    Now let $w\in\oC_S$ be weakly even. Here $\mydivisor{3}{-}{w}$ is
    big and nef, so we find the following table.
    \begin{equation*}
        \begin{array}{|c||c|c|c|}\hline
            \ww=6               & h^0 & h^1 & h^2 \\\hline\hline
            \mydivisor{ }{-}{w} &  1  &  0  &  0  \\\hline
            \mydivisor{3}{-}{w} &  5  &  0  &  0  \\\hline
        \end{array}
        \quad\quad   
        \begin{array}{|c||c|c|c|}\hline
            \ww=10              & h^0 & h^1 & h^2 \\\hline\hline
            \mydivisor{ }{-}{w} &  0  &  0  &  0  \\\hline
            \mydivisor{3}{-}{w} &  4  &  0  &  0  \\\hline
        \end{array}
    \end{equation*}
\end{example}
\section{Concluding remarks}
It is very likely that theorem \ref{theorem:main} is true for surfaces
of arbitrary degree. Unfortunately I cannot prove this. The main
obstruction is to exclude the
possibility that an irreducible contact surface is singular along
a curve which is contained in the nodal surface.

In \cite{barth2} Barth gave a construction of nodal surfaces
admitting even sets of nodes. The surfaces are constructed
as degeneracy locus of a generic quadratic form on a
globally generated vector bundle
on $\PP_3$.
For convenience we give a list of the strictly even sets of nodes
which have been obtained so far. Note that Barth's construction
gives exactly the gap of corollary \ref{corollary:main}.
\begin{center}
    \begin{tabular}{l|l}
        degree & even sets \\\hline
        3      & 4         \\
        4      & 8,16      \\
        5      & 16,20     \\
        6      & 24,32,40  \\
        8      & 48,64,72,80,\ldots,128 \\
        10     & 80,120,128,136,\ldots,208
    \end{tabular}
\end{center}
\noindent%
Stephan Endra\ss\\
Johannes Gutenberg-Universit\"at\\
Fachbereich 17\\
Staudinger Weg 9\\
D-55099 Mainz\\
{\tt endrass@mathematik.uni-mainz.de}
\end{document}